\documentclass[]{article}
\usepackage{graphicx}
\usepackage{graphics}
\usepackage{amssymb,amsmath,amsthm,epsfig}
\newcommand{\J}{P}
\newcommand{\TT}{\mathbb T}
\newcommand{\wJ}{\widetilde P}

\newcommand\dint{\displaystyle\int}
\newcommand\ds{\displaystyle}
\newtheorem{theorem}{Theorem}
\newtheorem{lemma}{Lemma}

\newtheorem{proposition}{Proposition}

\newtheorem{remark}{Remark}
\newcommand{\norm}[1]{{\left\|{#1}\right\|}}

\newcounter{reh}
\newcounter{rek}
\begin{document}

\begin{center}
{\large {\bf  Discrete Prolate Spheroidal Wave Functions: Further spectral analysis and some related applications.}}\\
\vskip 1cm Mourad Boulsane$^a$,  NourElHouda Bourguiba$^a$   and Abderrazek Karoui$^a$ {\footnote{
Corresponding author: Abderrazek Karoui, Email: abderrazek.karoui@fsb.rnu.tn\\
This work was supported by the Tunisian DGRST  research grant  UR 13ES47.}}
\end{center}
\vskip 0.25cm {\noindent
$^a$  University of Carthage, Department of Mathematics, Faculty of Sciences of Bizerte, Jarzouna 7021,  Tunisia.}\\

\noindent{\bf Abstract}--- For fixed $W\in \big(0,\frac{1}{2}\big)$ and positive integer $N\geq 1,$ the discrete   prolate spheroidal wave
functions (DPSWFs), denoted by  $U_{k,W}^N,$ $0\leq k\leq N-1$ form the set of  the eigenfunctions of the positive and  finite rank integral operator 
$\widetilde Q_{N,W},$ defined 
on $L^2(-1/2,1/2),$ with kernel $K_N(x,y)=\frac{\sin(N\pi(x-y))}{\sin(\pi(x-y))}\, \mathbf 1_{[-W,W]}(y).$  It is well known that the DPSWF's have a wide range of classical as well as recent signal processing applications. These applications
rely heavily on the properties of the DPSWFs as well as the behaviour of their eigenvalues  $\widetilde \lambda_{k,N}(W).$  In his pioneer work \cite{Slepian}, D. Slepian has given  the properties of the DPSWFs, their asymptotic approximations as well as the asymptotic behaviour and asymptotic decay rate of these eigenvalues. In this work, we give further properties as well as new non-asymptotic decay rates of the spectrum of the operator $\widetilde Q_{N,W}.$ In particular, we show that each eigenvalue $\widetilde \lambda_{k,N}(W)$ is up to a small constant  bounded above by the corresponding eigenvalue, associated with the classical prolate spheroidal wave functions (PSWFs). Then, based on the well established 
 results concerning the distribution and the decay rates of the eigenvalues associated with the PSWFs, we extend these results to 
the eigenvalues $\widetilde \lambda_{k,N}(W)$. Also, we show that the DPSWFs can be used for the  approximation of classical band-limited functions and they are well adapted for the approximation of functions from periodic Sobolev spaces. Finally, we provide the reader with some numerical examples that illustrate the different results of this work.\\



\section{Introduction}

A breakthrough in the theory and the construction of the discrete prolate spheroidal wave functions  is due to D. Slepian \cite{Slepian}, who has studied most of the properties,  the numerical computations, as well as the asymptotic behaviours of the DPSWFs and their  associated eigenvalues. Note that for fixed $W\in \big(0,\frac{1}{2}\big)$ and   integers $N_0\in \mathbb{N}, N\geq 1$, the DPSWF's are characterized as the amplitude spectra (Fourier series) of  index-limited complex  sequences with index support $[[N_0, N_0+N-1]]=\{N_0,\ldots,N_0+N-1\},$ that are most concentrated in the interval $(-W,W).$ 
For the sake of simplicity of the notations an without loss of generality, we  will only consider the  $N_0=0$  in this work. As it will described later on, the DPSWFs's are closely related to their associated Discrete Prolate Spheroidal Sequences (DPSS's). These DPSS's are infinite  sequences in $\ell_2(\mathbb C)$ with amplitude spectra supported in $[-W,W]$ and with  coefficients most concentrated in the index range $[[0,\ldots,N-1]].$
The DPSS's sequences have been successfully used in various classical as well as fairly  recent applications from the signal processing area.
To cite but a few, the prediction of white noise random samples of  discrete signals with bandwidth $W_0,$ \cite{Slepian},  the DPSS's based scheme for compressive sensing \cite{Wakin1},  parametric waveform and detection of extended targets \cite{Feng} and fast algorithms for Fourier extension \cite{Adcock}, etc.

It has been shown in \cite{Slepian}, that the solution of the energy maximization problem, associated with the DPSWF's  is given by the first eigenfunction corresponding to the largest eigenvalue of the following eigenproblem 
\begin{eqnarray}
\ds \int_{-W}^{W}
\frac{\sin(N\pi(y-x))}{\sin(\pi(y-x))} {h}(y)dy&=&\lambda {h}(x),\quad x\in (-1/2,1/2).
\end{eqnarray}
Therefore,  the different DPSWF's $U_{k,W}^N$ are the $N$ eigenfunctions of a finite rank integral operator $\widetilde{Q}_{W,N},$ that is  
\begin{equation}\label{eigenproblem1}
\widetilde{Q}_{W,N}(U_{k,W}^N)(x)=\int_{-W}^{W}\frac{\sin(N\pi(x-y))}{\sin(\pi(x-y))}U_{k,W}^N(y)dy=\widetilde{\lambda}_{k,N}(W)U_{k,W}^N(x).
\end{equation}
Here, $1>\widetilde{\lambda}_{0,N}(W)>\widetilde{\lambda}_{1,N}(W)>\cdots >\widetilde{\lambda}_{N-1,N}(W)$ is the sequence of the associated eigenvalues, arranged in the decreasing order. The $N$ DPSWF's  form an orthonormal system of  both  $L^2(-W,W),\, W\in \big(0,\frac{1}{2}\big)$ and $L^2(-1/2,1/2).$ More precisely, they satisfy  the following double orthogonality properties 
\begin{equation}
\label{doubleorthogonality}
\int_{-W}^{W} U_{k,W}^N(x) \overline{U_{j,W}^N(x)}dx=\widetilde{\lambda}_{k,N}(W)\delta_{k,j},\quad \int_{-1/2}^{1/2} U_{k,W}^N(x) \overline{U_{j,W}^N(x)}dx=\delta_{k,j},\quad k, j = 0,\ldots,N-1.
\end{equation}
From \cite{Slepian}, the DPSWFs are related to the DPSS's by the following rule. Let  $V_{k,W}^N=(v_0^{(k)},.....,v_{N-1}^{(k)})^{T}, $ $k=0,..,N-1$ be the $N$ vectors obtained by truncating the DPSS's to the index set $[[0,N-1]].$ Then, these truncated DPSS's are the $N$ eigenvectors of the Toeplitz  matrix
\begin{equation}\label{Toeplitzmatrix}
\rho_{N,W}=\left[\frac{\sin(2\pi(n-m) W)}{\pi(n-m)}\right]_{n,m=0,..,N-1}.
\end{equation}
Moreover, we have 
\begin{eqnarray}\label{DPSWFs}
U_{k,W}^N(x) &=& \epsilon_k \ds{\sum_{n=0}^{N-1}} v_n^{(k)}e^{-i \pi(N-1-2n)x},\qquad \epsilon_k=\left\{\begin{array}{ll}
1, & \hbox{k even;} \\
i, & \hbox{k odd.}
\end{array}
\right.
\end{eqnarray}
Note that the matrix $\rho_{N,W}$ has the same spectrum as the integral operator $\widetilde Q_{W,N},$ that is the DPSWFs  $U_{k,W}^N$
and the corresponding truncated  DPSS's $V_{k,W}^N$ are associated with same eigenvalues $\widetilde \lambda_{k,N}(W).$
 Also, it is interesting to note that the spectrum associated with the DPSWFs has some surprising similarities with the spectrum associated with the classical PSWFs, that have been introduced and greatly investigated since the early 1960's, by D. Slepian and his co-authors H. Landau and H. Pollak, see \cite{Landau, Sl1, Slepian}. We recall that for a given real number  $c>0$, called the bandwidth, the  PSWFs $(\psi_{n,c}(\cdot))_{n\geq 0}$ constitute an orthonormal basis of $L^2([-1, +1]),$ an orthogonal system of $L^2({\bf R})$ and an
orthogonal  basis of  the Paley-Wiener space $B_c,$ given by 
\begin{equation}\label{Bc}
B_c=\left\{ f\in L^2({\mathbb R}),\,\, \mbox{Support\ \ } \widehat f\subset [-c,c]\right\}.
\end{equation}
Here, $\widehat f$ denotes the Fourier transform of $f\in L^2(\mathbb R).$ They are eigenfunctions of the
 Sinc-kernel operator  defined  on $L^2([-1,1]),$ that is 
\begin{equation}\label{eq1.1}
\mathcal Q_c(\psi_{n,c})(x)=\int_{-1}^1\frac{\sin c(x-y)}{\pi (x-y)}\, \psi_{n,c}(y)\, dy= \lambda_n(c) \psi_{n,c}(x),\quad x\in [-1,1].
\end{equation}
Unlike the classical case, where there exist a rich literature on the behaviour and the decay rates (both asymptotic and non-asymptotic) 
of the eigenvalues $\lambda_n(c),$ see for example \cite{Bonami-Karoui1, BJK, HL, Landau, Sl1, Slepian2}, the counterpart  literature 
for the $\widetilde \lambda_{k,N}(W)$ is still very limited. The main existing decay rate result for  the $\widetilde \lambda_{k,N}(W)$ is an asymptotic
one and it  goes back to \cite{Slepian}, where it has been shown that for fixed $W\in \big(0,\frac{1}{2}\big)$ and $\varepsilon \in \big(0,\frac{1}{2W}-1),$ we have 
\begin{equation}\label{decay1}
\widetilde \lambda_{k,N}(W) \leq C_1(W,\varepsilon) e^{-C_2(W,\varepsilon) N},\quad \forall\, k\geq \lceil{2NW(1+\varepsilon)}\rceil,\quad N \geq N_1(W,\varepsilon),
\end{equation}
for some constants  $C_1(W,\varepsilon),$ $C_2(W,\varepsilon)$ and $N_1(W,\varepsilon)\in \mathbb N$ that depend on $W, \varepsilon.$  
The previous estimate is asymptotic and the dependence of the previous constants does not have explicit estimates. 
The previous decay rate has been recently generalized in \cite{zhu} to the multiband DPSS's setting. Moreover, in this last reference and by using some 
advanced matrix analysis and computations techniques, the authors have given the following distribution of the $\widetilde \lambda_{k,N}(W).$
If $\mathbb W$ is an union of $J$ pairwise disjoint intervals with $\mathbb W\subset \big(-\frac{1}{2},\frac{1}{2}\big)$ and $\varepsilon \in 
\big(0,\frac{1}{2}\big),$ then
\begin{equation}\label{decay2}
\#\{k: \varepsilon \leq \widetilde \lambda_{k,N}(\mathbb W)\leq 1-\varepsilon \} \leq J \frac{\frac{2}{\pi^2} \log(N-1)+\frac{2}{\pi^2} \frac{2N-1}{N-1}}{\varepsilon (1-\varepsilon)}.
\end{equation}
In this work, for $c= \pi N W$ and by comparing the Hilbert-Schmidt norms 
$\|\widetilde Q_{W,N}\|_{{\scriptscriptstyle \rm HS}}$ and $\|\mathcal  Q_{c}\|_{{\scriptscriptstyle \rm HS}}$ of the integral operators $\widetilde Q_{W,N}$ and $\mathcal Q_c,$ given by \eqref{eigenproblem1} and \eqref{eq1.1}, we prove that for $J=1,$ we have 
\begin{equation}\label{decay3}
\#\{k: \varepsilon \leq \widetilde \lambda_{k,N}(W)\leq 1-\varepsilon \} \leq \frac{\frac{1}{\pi^2} \log(2 N W)+0.45-   \frac{2}{3}W^2 +\frac{W^2}{6c^2}\sin^2(2c)}{\varepsilon(1-\varepsilon)},\qquad c=\pi N W.
\end{equation}
It can be easily checked that for $J=1,$ $N\geq 2$ and $\pi N W\geq 1,$ our estimate \eqref{decay3} improves the estimate \eqref{decay2}.
Also, the comparison of the previous Hilbert-Schmidt norms, together with the use of the Wielandt-Hoffman inequality, we prove that for sufficiently small $W,$ the spectrum of $\widetilde Q_{W,N}$ is well  approximated in the $\ell_2$-norm by the spectrum of the Sinc-kernel operator $\mathcal Q_c,\, c=\pi NW.$ More precisely, for any $N\geq 1$ and $W\in \big(0,\frac{1}{2}\big),$ we have
\begin{equation}
\label{approx1}
\left(\ds  \sum_{k=0}^{N-1}\left|\widetilde{\lambda}_{k,N}(W)-\lambda_{k}(c)\right|^2\right)^{\frac{1}{2}} \leq  W^3\left( \frac{4\pi^2}
{3\sin(2 W\pi)}\right),\qquad c= \pi N W.
\end{equation}
Also by taking advantage from a connection between the energy maximization problems associated to the DPSWFs and the classical PSWFs, $\psi_{n,c}$ with $c= \pi N W,$ we prove the following unexpected and important result relating the $\widetilde\lambda_{n,W}^N$ and the $\lambda_n(c),$
\begin{equation}\label{decay5}
\widetilde{\lambda}_{n,N}(W) \leq  A_W \,\,\,\lambda_n(c), \quad 0\leq n\leq N-1,
\end{equation}
where
\begin{equation}\label{ddecay5}
 \frac{\pi^2}{8}\leq A_W = \frac{2 \pi^2}{\cos^2(\pi W)} \left(\frac{1}{4}-W^2\right)^2 \leq 2.
\end{equation}
Thanks to the estimate \eqref{decay5}, all the existing known asymptotic and non-asymptotic decay rates for the $\lambda_n(c)$ are transmitted 
to the $\widetilde{\lambda}_{n,N}(W).$ For example, based on the recent non-asymptotic estimates of the $\lambda_n(c),$ given in \cite{BJK}, one concludes that under the condition that for sufficiently away from the plunge region of the spectrum, that is for $2\leq \frac{e \pi}{2} N W \leq n\leq N-1,$ we have 
\begin{equation}\label{decay6}
\widetilde{\lambda}_{n,N}(W) \leq  2 \,\,\, e^{-(2n+1) \log\big(\frac{2n+2}{e\pi NW}\big)}.
\end{equation}
Moreover, for $n$ close to the plunge region around $[2NW],$ there exists a constant $\eta >0,$ such that 
\begin{equation}\label{decay7}
\widetilde{\lambda}_{n,N}(W) \leq 2 e^{- \eta \frac{n-2NW}{\log(\pi NW)+5}},\qquad 2NW +\log(\pi N W)+6\leq n\leq \pi N W.
\end{equation}

As applications of the DPSWF's that  we consider in this work, we first get an estimate of the unknown constant appearing in the Tur\`an-Nazarov concentration inequality. Then, we check that there exists $N_0\geq [2NW]-1,$ such that  the eigen-space spanned by the first $N_0$ dilated DPSWFs $\sqrt{W} U_{k,W}^N(W\cdot)$ is approximated by the eigen-space spanned by the corresponding classical $\psi_{k,c}.$ Also, we check that these 
DPSWF's are well adapted for the spectral approximation of functions from the periodic Sobolev space $H^s_{per}(-1/2,1/2),\, s>0.$

Finally, this work is organized as follows. In section $2$, we give some mathematical preliminaries related to the properties and the computations of DPSWFs and DPSS's and their  associated eigenvalues.  Moreover, we give some first estimates for the eigenvalues  associated with the DPSWF's. In section 3, we study some interesting connections between the DPSWFs and their corresponding
classical $\psi_{n,c},\, c=\pi NW.$ Based on these connections, we deduce various results on the distribution and the decay rates of 
the eigenvalues $\widetilde \lambda_{n,N}(W).$ Section 4 is devoted to the previous proposed applications of the DPSWF's. In the last section 5, we give some numerical examples  that illustrate the different results of this work.

\section{Mathematical preliminaries}

In this paragraph, we first recall from the literature, some properties and computational methods for the DPSWFs and their associated eigenvalues.
Also, we give some first  estimates of the eigenvalues associated with the DPSWFs. These estimates are obtained in a fairly easy way by using the Min-Max characterization of the eigenvalues of self-adjoint compact operators.  More involved and precise estimates of the eigenvalues $\widetilde \lambda_{n,N}(W),$ is the subject of the next section 3.

We recall from \cite{Slepian},  that the DPSWFs $(U_{n,W}^N)_{0\leq n\leq N-1}$  are the eigenfunctions of the positive, self-adjoint  finite rank  integral operator $\widetilde{Q}_{W,N},$ given by \eqref{eigenproblem1}. This last eigen-problem is a consequence of the fact that the DPSWF's. Among 
the space $S_N$ of all sequences  ${\pmb x}=(x_n)_n \in l^2(\mathbb C)$ with elements indexed on 
$[[0,N]],$ so that their amplitude spectra  ${\displaystyle  \widehat{{\pmb x}}(t)= \sum_{k=0}^N x_k e^{2i\pi k t},}$ find those sequences
with amplitude spectra most concentrated on $(-W,W),$ that is solve the maximization problem
\begin{equation}\label{Maximization1}
U = \arg \max_{{\pmb x}\in S_N} \frac{\|\widehat{{\pmb x}}(t)\|^2_{L^2(-W,W)}}{\|\widehat {{\pmb x}}(t)\|^2_{L^2(-1/2,1/2)}}.
\end{equation}
Note that the DPSWFs $(U_{n,W}^N)$ are periodic. They have period 1 if $N$ is odd and period 2 if $N$ is even. In either case we have
$$U_{k,W}^N(x+1) = (-1)^{N-1}U_{k,W}^N(x),\qquad x\in [-1/2,1/2).$$ 
Also, the associated eigenvalues satisfy the following relation,
\begin{equation}
\widetilde \lambda_{k,N}\Big(\frac{1}{2}-W\Big)=1-\widetilde \lambda_{N-k-1,N}(W),\quad \forall\, k=0,\ldots,N-1.
\end{equation}

Moreover, the DPSWFs can be computed by using two schemes. The first scheme is given by \eqref{DPSWFs}, that is an expansion with respect 
to the eigenvectors of the Toeplitz matrix $\rho_{N,W},$ given by \eqref{Toeplitzmatrix}. Note that from \cite{Slepian} and  \cite{zhu}, the  matrix
$\rho_{N,W}$ is the matrix representation of $\mathcal I_N \mathcal B_W \mathcal I_n^*,$  a composition of  index- and band-limiting operators,  $\mathcal I_N: \ell_2(\mathbb C)\rightarrow \mathbb C^N,\,\,  \mathbb B_W:\ell_2(\mathbb C)\rightarrow \ell_2(\mathbb C) ,$ given for $h=(h_n)_{n\in \mathbb Z},$ by 
$$  \mathbb B_W (h)(m)= \sum_{n\in \mathbb Z} \frac{\sin (2\pi W (m-n))}{\pi (m-n)} h_n,\qquad \mathcal I_N (h)(m)= h_m,\quad m\in \{0,\ldots, N-1\}.$$
This is a consequence of the connection between the DPSWF's and their associated DPSS's.  
We should mention that the DPSS's are solutions of the following energy maximization
dual problem. Let  $\mathcal{B}_{W}$ be the Paley-Wiener space given by $\mathcal{B}_{W}=\{ h=(h_n)_{n\in \mathbb{Z}}\in \ell^2(\mathbb{C}),\,\, \mbox{Supp }(\widehat{h})\subset [-W,W] \}.$ Here, $\widehat{h}(x)=\ds \sum_{n\in\mathbb{Z}} h_n e^{-i\pi (N-1-2n)x}, x\in (-1/2,1/2).$ Then, find 
$$  h =\ds \arg\max_{h\in \mathcal{B}_{W}} \sum_{n=0}^{N-1}|h_n|^2/{\Big(\sum_{n=-\infty}^{+\infty}|h_n|^2\Big)}.$$
Consequently, the  DPSS's are solutions of the system of equations
\begin{equation}\label{DPSSs}
 \sum_{m\in \mathbb Z} \frac{\sin (2\pi W (m-n))}{\pi (m-n)} v_m^{(k)}= \widetilde \lambda_{k,N}(W) v_n^{(k)},\quad \forall n\in \mathbb Z.
\end{equation}
For more details, see \cite{Slepian}.

The second scheme for the construction of the DPSWF's is  based on the computation of the eigenvectors of a Sturm-Liouville differential operator $M_{W,N},$ 
commuting with  the integral operator $\widetilde Q_{W,N},$ see for example \cite{Slepian}. This  differential operator is given by 
\begin{equation}\label{differentialoperator}
    M_{W,N}(g)(x)=\frac{1}{4\pi^2}\frac{d}{dx}\left[ \left(\cos(2\pi x)-A\right)\frac{d}{dx}(g)(x)\right]+\frac{1}{4}(N^2-1)\cos(2\pi x)(g)(x),\quad 
    A=\cos(2\pi W).
\end{equation}
Hence, the DPSWFs are also given in terms with the eigenvectors of $M_{W,N}.$  It is easy to check that 
\begin{eqnarray}\label{differentialoperator2}
 M_{W,N}(e^{i\pi (N-1-2n)x})(x) &=&\frac{1}{2}n(N-n)e^{i\pi (N-2n+1)x} +
 \left[A \left(\frac{N-1}{2}-n\right)^2\right]e^{i\pi (N-2n-1)x}\nonumber\\
&+&\frac{1}{2}(n+1)(N-n-1)e^{i\pi (N-2n-3)x}\label{3}
\end{eqnarray}
Consequently, the expansion coefficients in the basis $\{e^{i\pi (N-1-2k)x},\,\,\, 0\leq k\leq N-1\}$ of the  $n$-th  DPSWFs $(U_{k,W}^N)$ 
are given by the components of  the $n$-th  eigenvector of the $N\times N$ tri-diadiagonal matrix $\sigma(N,W),$ with coefficients given by $$\sigma(N,W)_{ij}=\left\{
  \begin{array}{ll}
    \frac{1}{2}i(N-i), & \hbox{$j=i-1;$} \\
    \cos(2\pi W) \left(\frac{N-1}{2}-i\right)^2 & \hbox{$j=i;$} \\
    \frac{1}{2}(i+1)(N-i-1), & \hbox{$j=i+1;$} \\
    0, & \hbox{ $|j-i| >1,$}
  \end{array}
\right.\qquad i,j =0,\ldots,N-1.$$

It is interesting to note that by considering the finite rank and  positive-definite integral operator $\widetilde Q_{W,N}$ as an operator acting on the 
Hilbert space $L^2(-W,W)$ and by using the Min-Max theorem for this operator, one gets the following lemma that provides us with a partial result related to the decay rate of the $\widetilde Q_{W,N}.$ We should mention that the proof of this lemma  mimics the technique used in \cite{BJK} for proving a similar result concerning a  decay rate of the $\lambda_n(c),$ the eigenvalues of the operator $\mathcal Q_c,$ given by \eqref{eq1.1} and associated with the classical PSWFs.

\begin{lemma}
 For any real number $0<W<\frac{2}{e \pi}$ and any integer $N\geq 2,$   we have
\begin{equation}\label{Eq3.1}
    \widetilde{\lambda}_{n,N}(W)\leq \frac{C_W}{\sqrt{N-1}\log\left(\frac{2n}{e\pi W(N-1)}\right)}
\left(\frac{e\pi W(N-1)}{2n}\right)^{n-\frac{1}{2}},\qquad  \frac{e\pi W(N-1)}{2} < n \leq N-1,
\end{equation}
where $C_W= \sqrt{2W} \Big(2+\frac{2}{e\pi W}\Big).$
\end{lemma}

\noindent
{\bf Proof:} We first recall the Courant-Fischer-Weyl Min-Max variational principle concerning the positive eigenvalues
of a self-adjoint compact operator $A$ acting on a Hilbert space $\mathcal H,$ with eigenvalues arranged in the decreasing order $\lambda_0\geq \lambda_1\geq \cdots\geq \lambda_n\geq \cdots . $ In this case, we have
$$\lambda_n = \min_{f\in S_n}\,\,\,  \max_{f\in S_n^{\perp},\|f\|_{\mathcal H}=1} < Af, f>_{\mathcal H},$$
 where $S_n$ is a subspace of $\mathcal H$ of dimension $n.$ In our case, we have $A=\widetilde{Q}_{W,N},$ $\mathcal H= L^2(-W,W).$ We consider the special case of  
 $$S_n=\mbox{Span}\left\{\wJ_0\left(x\right), \wJ_1\left(x\right),\ldots,\wJ_{n-1}\left(x\right)\right\}$$
  and
 $$f(x)=\ds\sum_{k\geq n} a_k \wJ_k \left(x\right)\in S_n^{\perp},\qquad
  \parallel f\parallel_{L^2_{([-W,W])}}=\ds\sum_{k\geq n}|a_k|^2=1.$$
  Here, ${\displaystyle \wJ_k(x)= \sqrt{\frac{2k+1}{2W}}\J_k \left(\frac{x}{W}\right)},$ where $\J_k$ is the usual Legendre polynomial of degree $k$ and satisfying $P_k(1)=1.$ Note that the $\wJ_k$ form an orthonormal family of  $L^2_{(-W,W)}.$ The normalization constant follows from the fact that  
\begin{eqnarray}\label{Legendre}
  \ds{\parallel \J_k(\cdot)\parallel_{L^2_{(-W,W)}}} &=& \ds{\left(\ds\int_{-W}^{W }| \J_k \left(\frac{y}{W}\right)|^2dy\right)^{\frac{1}{2}}= \left(W\ds\int_{-1}^{1 }| \J_k (y)|^2dy\right)^{\frac{1}{2}}=\sqrt{\frac{2W}{2k+1}}=h_{k,W}}.
\end{eqnarray}
On the other hand, we have 
\begin{eqnarray}\label{Legendre2}
  \widetilde{Q}_{W,N} \left(\wJ_k \right)(x) &=&\ds\int_{-W}^{W}\frac{\sin(N\pi(x-y))}{\sin(\pi(x-y))} \wJ_k \left(y\right)dy =\ds\int_{-W}^{W}\ds\sum_{j=0}^{N-1}\frac{e^{i\pi(N-1-2j)(x-y)}}{h_{k,W}} \J_k \left(\frac{y}{W}\right)dy  \nonumber \\
   &=&W \ds\sum_{j=0}^{N-1}e^{i\pi(N-1-2j)x}\ds\int_{-1}^{1}\frac{e^{-iW\pi(N-1-2j)y}}{h_{k,W}} \J_k (y)dy
\end{eqnarray}
Moreover,  it is known that, see for example \cite{NIST}
 \begin{eqnarray}\label{Eq2.1}
    \dint_{-1}^1 e^{ixy}\J_k(y) d y&=&i^k\sqrt{\frac{2\pi}{x}} J_{k+\frac{1}{2}} (x)
   ,\quad x\in \mathbb R.
  \end{eqnarray}
  where $J_{\alpha}$ is the Bessel function of the first type and order $\alpha > −1.$
Further, the Bessel function $J_{\alpha}$ has the following fast decay with respect to the parameter $\alpha$,
  \begin{eqnarray}
    |J_{\alpha}(z)| &\leq& \frac{\big|\frac{z}{2}\big|^\alpha}{\Gamma(\alpha+1)}
  \end{eqnarray}
 Here, $\Gamma(\cdot)$ is the Gamma function, that satisfies the following bounds, see \cite{Batir} that
\begin{equation}\label{Eq2.2}
    \sqrt{2e}\left(\frac{x+\frac{1}{2}}{e}\right)^{x+\frac{1}{2}}\leq \Gamma(x+1) \leq  \sqrt{2\pi}\left(\frac{x+\frac{1}{2}}{e}\right)^{x+\frac{1}{2}},\quad x>-1.
\end{equation}
From the previous inequality and \eqref{Legendre}, we deduce that
\begin{eqnarray}\label{Legendre3}
  \Big|\ds\int_{-1}^{1}\frac{e^{-iW\pi(N-1-2j)y}}{h_{k,W}} \J_k (y)dy  \Big| &\leq& \sqrt{\frac{2k+1}{W^2|N-1-2j|2(k+1)}}\left(\frac{e W \pi |N-1-2j|}{2(k+1)}\right)^{k+\frac{1}{2}}
\end{eqnarray}
Then, by using \eqref{Legendre2}, \eqref{Legendre3} and  the Minkowski's inequality, one gets for $k\geq e\pi W (N-1)/2,$
\begin{eqnarray}
 \left\|\widetilde{Q}_{W,N}\left(\wJ_k \right) \left(x\right)\right\|_{L^2(-W,W)}
   &\leq &\ds\sum_{j=0}^{N-1}\|{W e^{i\pi(N-1-2j)x}}\|_{L^2(-W,W)}\sqrt{\frac{2k+1}{W^2|N-1-2j|2(k+1)}}\left(\frac{e W \pi |N-1-2j|}{2(k+1)}\right)^{k+\frac{1}{2}}
   \nonumber\\
&\leq &\sqrt{2W} \ds\sum_{j=0}^{N-1}
   \sqrt{\frac{e\pi W}{2(k+1)}}\left(\frac{e W \pi |N-1-2j|}{2(k+1)}\right)^{k} \nonumber\\
&\leq &\frac{C_W}{\sqrt{N-1}}\left(\frac{e W \pi (N-1)}{2(k+1)}\right)^{k+\frac{1}{2}},\quad C_W=\sqrt{2W} \Big(2+\frac{2}{e\pi W}\Big).
\end{eqnarray}
The last inequality follows from the fact that for  $k\geq e\pi W (N-1)/2,$
$$\sum_{j=0}^{N-1} |N-1-2j|^k\leq 2(N-1)^k+\frac{(N-1)^{k+1}}{k+1}\leq (N-1)^k \Big(2+\frac{2}{e\pi W}\Big).$$
Hence, for the previous $f\in S_n^{\perp},$  and by using H\"older's inequality, and taking into account that  $\|f\|_{L^2(I,\omega_{W})}= 1,$ so that $|a_k|\leq 1,$ for $k\geq n,$  one gets
\begin{eqnarray}\label{Eq3.3}
|< \widetilde{Q}_{W,N}f,f>_{{L^2([-W,W])}}| &\leq & \sum_{k\geq n} |a_k| \|\widetilde{Q}_{W,N}
 \wJ_k\left(\cdot\right)\|_{L^2([-W,W])}\nonumber\\
&\leq &\frac{C_W}{\sqrt{N-1}}\sum_{k\geq n} |a_k| \left(\frac{eW \pi(N-1)}{2(k+1)}\right)^{k+\frac{1}{2}}
\leq \frac{C_W}{\sqrt{N-1}}\sum_{k\geq n}  \left(\frac{W e \pi(N-1)}{2(k+1)}\right)^{k+\frac{1}{2}}.
\end{eqnarray}
The decay of the sequence appearing in the previous sum, allows us to compare this later with its integral counterpart, that is
\begin{equation}\label{Eq3.4}
\sum_{k\geq n}  \left(\frac{W e \pi(N-1)}{2(k+1)}\right)^{k+\frac{1}{2}}\leq \int_{n-1}^{+\infty} e^{- (x+\frac{1}{2})\log\left(\frac{2(x+1)}{e W\pi (N-1)}\right)}\; dx\leq
\int_{n-1}^{+\infty} e^{- (x+\frac{1}{2})\log\left(\frac{2n}{e W\pi (N-1)}\right)}\, dx
\end{equation}
Hence, by using \eqref{Eq3.3} and \eqref{Eq3.4}, one concludes that
\begin{equation}\label{Eq3.5}
\max_{f\in S_n^{\perp},\, \|f\|_{L^2(I,\omega_{W})=1}} <\widetilde{Q}_{W,N} f,f>_{L^2([-W,W])}\leq
 \frac{C_W}{\sqrt{N-1}\log(\frac{2n}{e W\pi (N-1)})}
 e^{-(n-\frac{1}{2}) \log(\frac{2n}{e W\pi (N-1)})}.
\end{equation}
 To conclude for the proof of the lemma, it suffices to use the previous  Courant-Fischer-Weyl Min-Max variational principle.$\qquad \Box $

\section{The Spectrum associated with the DPSWF's: Behaviour and decay rates}

In the first part of this  section, we estimate the Hilbert-Schmidt norms of the two operators \eqref{eigenproblem1} and \eqref{eq1.1}. As consequences,
we  give a comparison in the $\ell_2$-norm of the  spectrum associated with the DPSWFs with parameters $N, W$ and the 
spectrum associated with the classical PSWFs $(\psi_{n,c})_n$ wih $c=N\pi W.$ Also, we give a fairly precise estimate of the number of the eigenvalues 
$\widetilde \lambda_{k,N}(W)$ lying in the interval $[\varepsilon, 1-\varepsilon],$ where $\varepsilon \in (0,1/2).$ 
In the second part, we use the energy maximizations characterizations of the DPSWFs and the PSWFs, and get an interesting fairly precise upper bound of   the eigenvalues $\widetilde \lambda_{n,N}(W)$ in terms of the eigenvalues $\lambda_n(c),$ for $0\leq n\leq N-1.$ As a consequence and by using 
the well established decay rates and behaviour of the $\lambda_n(c),$ we deduce similar results for the $\widetilde \lambda_{n,N}(W),$ with $c=N\pi W.$ The following proposition provides us with an $\ell_2$-estimate of the spectrum of $\widetilde Q_{W,N}$ by the spectrum of $\mathcal Q_c.$ Note that since 
$\widetilde Q_{W,N}$ is of finite rank, which is not the case for the operator $\mathcal Q_c,$ then this  $\ell_2$-estimate is done under the rule that 
$\widetilde \lambda_{k,N}(W)=0,$ whenever $k\geq N.$

\begin{proposition}\label{OperatorsHSnorms}
Under the previous notation, for $W \in (0,\frac{1}{2})$ and an integer  $N\geq 1,$ we have for $c= \pi N W$
\begin{equation}\label{6}
\| \lambda(\widetilde Q_{W,N})-\lambda (\mathcal Q_c)\|_{\ell_2}=\left(\ds  \sum_{k=0}^{\infty}\left|\widetilde{\lambda}_{k,N}(W)-\lambda_{k}(c)\right|^2\right)^{\frac{1}{2}} \leq  W^3\left( \frac{4\pi^2}
{3\sin(2 W\pi)}\right)
\end{equation}
\end{proposition}

\noindent
{\bf Proof:}
Since the operator $\widetilde{Q}_{W,N}$ acts on $L^2(-W,W)$. Then we consider the operator $\mathcal{Q}_{W,c}$ associated with the classical
PSWFs that are mostly concentrated on $[-W,W] $ and have  bandwidth $[-c,c] $. These last family of PSWFs are solutions of the eigenvalues problem
\begin{equation}
\label{Q1}
\mathcal Q_{W,c}(\psi)=\ds\int_{-W}^{W}\frac{\sin(c(x-y))}{\pi(x-y)}\psi_{k,W}(y)dy=\lambda_{k,W}(c)\psi_{k,W}(x).
\end{equation}
It is well know that $\lambda_{k,W}(c)=\lambda_{k,1}(cW)=\lambda_k(c W),$ $\forall W>0.$
It is common to write $\lambda_{k,1}(c)=\lambda_{k}(c)$ and $\mathcal Q_{1,c}=\mathcal Q_c,$ where this later is given by \eqref{eq1.1}.
For $W\in (0,\frac{1}{2})$, we let  $c_N= \pi N $ and $c= \pi N  W.$ Then, we have
\begin{eqnarray}\label{Eq3}
  \|\widetilde{Q}_{W,N}-Q_{W,c_N}\|^2_{HS} &=& \ds\int_{-W}^{W}\ds\int_{-W}^{W}\left(\frac{\sin(c_N(x-y))} {\sin(\pi(x-y))}
-\frac{\sin(c_N(x-y))} {(\pi(x-y))}    \right)^2 dx  dy\nonumber\\
&=&W^2\ds\int_{-1}^{1}\ds\int_{-1}^{1} \left(\sin(c(t-u))\right)^2\left[ \frac{1} {\sin(W\pi(t-u))}
-\frac{1} {(W\pi(t-u))} \right]^2 du dt
\end{eqnarray}
But for $X=W\pi (t-u)\in[-2 \pi W,2 \pi W],$ we have
\begin{eqnarray}\label{Eq4}
  \left| \frac{X-\sin X}{X\sin X}\right| &\leq & \frac{|X|}{6}\left|\frac{X}{ \sin X} \right|\leq\frac{W\pi}{3}\frac{2 W \pi}{\sin(2W\pi)}.
\end{eqnarray}
The last inequality is due to the fact that $x\mapsto\frac{x}{\sin x}$ is increasing on $[-2\pi W,2\pi W].$ Consequently, by using 
\eqref{Eq3} and \eqref{Eq4}, one gets 
 \begin{eqnarray}\label{normsHS}
  \|\widetilde{Q}_{W,N}-Q_{W,c_N}\|^2_{HS} &\leq &  4 W^2\left( W^2\frac{2\pi^2}{3\sin(2\pi W )} \right)^2.
\end{eqnarray}
Finally, by using \eqref{Q1} and the previous equality together with  Wielandt-Hoffman inequality, one gets
\begin{eqnarray}\label{HSnorms}
\| \lambda(\widetilde Q_{W,N})-\lambda (\mathcal Q_c)\|_{\ell_2}=\left(\ds  \sum_{k=0}^{\infty}\left|\widetilde{\lambda}_{k,N}(W)-\lambda_{k}(c)\right|^2\right)^{\frac{1}{2}} &\leq &
W^3\left( \frac{4\pi^2}
{3\sin(2 W\pi)}\right).\quad \Box
\end{eqnarray}

Next by comparing the Hilbert-Schmidt norms of the operators $\widetilde Q_{W,N}$ and $\mathcal Q_c,$  together with a precise  estimate
of $ Trace(\mathcal Q_c)- \|{Q}_{c}\|^2_{HS},$ we get the following theorem, showing that the eigenvalues $\widetilde \lambda_{k,W}$ cluster
around $1$ and $0.$
 
\begin{theorem}\label{numbereigenvales}
For any $\varepsilon \in(0,1/2)$ and any $W\in(0,\frac{1}{2})$, let
$$ \mathcal N(W,\varepsilon) = \# \{k\; ; \; \varepsilon<\widetilde{\lambda}_{k,N}(W)<1-\varepsilon\},$$ then we have
\begin{eqnarray}\label{Eigenvaluesnumber}
\mathcal N(W,\varepsilon)   &\leq&\frac{\frac{1}{\pi^2} \log(2 N W)+0.45-   \frac{2}{3}W^2 +\frac{W^2}{6c^2}\sin^2(2c)}{\varepsilon(1-\varepsilon)},\qquad c=\pi N W.
\end{eqnarray}
\end{theorem}

\noindent
{\bf Proof:}  Since
 $ \|\widetilde{Q}_{W,N}\|^2_{HS} = \ds\int_{-W}^{W}\ds\int_{-W}^{W}\left(\frac{\sin(N\pi (x-y))} {\sin(\pi(x-y))}
    \right)^2dxdy$, then using the new variables $t=\frac{x}{W}, \; \; u=\frac{y}{W}$, we get
$$
   \|\widetilde{Q}_{W,N}\|^2_{HS} = W^2\ds\int_{-1}^{1}\ds\int_{-1}^{1}\left(\frac{\sin(\pi N W(u-x))} {\sin(\pi W(u-x))}
    \right)^2 du dx.
$$
That is for  $c= \pi N  W,$ we have
\begin{eqnarray}
  \|\widetilde{Q}_{W,N}\|^2_{HS}- \|\mathcal {Q}_{c}\|^2_{HS} &=&\ds\int_{-1}^{1}\left(\ds\int_{-1-x}^{1-x}(\sin(ct))^2\left(\frac{W^2} {\sin^2(\pi W t)}-\frac{1}{t^2\pi^2}
    \right)dt\right) dx\nonumber\\
&=&\ds\int_{-1}^{1}\left(\ds\int_{-1-x}^{1-x}(\sin(ct))^2h_{W}(t)dt\right) dx \label{8}
\end{eqnarray}
with $h_{W}(t)=W^2\left(\frac{1} {\sin^2(y)}-\frac{1}{y^2}
    \right)=W^2g(y)$ and $y \in  ]-\pi,\pi[$. We check that $h_{W}(t)\geq \frac{W^2}{3}$. In fact,   $g$ is even and increasing function on $[0,\pi].$
Note that straightforward computation gives us
${  g'(y) = \frac{2}{(\sin(y))^3}\left[ -\cos(y)+\left(\frac{\sin(y)}{y}\right)^3\right] }.$
It is clear that $g'(y)\geq0 $ if $\;y\in[\frac{\pi}{2},\pi],$ and  
\begin{eqnarray*}
   \left(\frac{\sin(y)}{y}\right)^3-\cos(y) &\geq& \Big(1-\frac{y^2}{6}\Big)^3 -\Big(1-\frac{y^2}{2}+\frac{y^4}{24}\Big)\geq \frac{y^4}{24}\Big(1-\frac{y^2}{9}\Big)\geq0,\quad 
  y\in \left[0,\frac{\pi}{2}\right]. 
\end{eqnarray*}
Consequently, we have 
\begin{equation}\label{7}
g(y)\geq \ds \inf_{y\in[0,\pi]}g(y)=\ds\lim_{y\rightarrow 0}g(y)=\frac{1}{3}.
\end{equation}
By combining \eqref{8} and \eqref{7}, one gets
\begin{eqnarray*}
  \|\widetilde{Q}_{W,N}\|^2_{HS}- \|\mathcal {Q}_{c}\|^2_{HS} &\geq &\frac{W^2}{3}
\ds\int_{-1}^{1}\left(\ds\int_{-1-x}^{1-x}\frac{1-\cos(2ct)}{2}dt\right) dx =
\frac{2 W^2}{3}-\frac{ W^2}{12 c}
\ds\int_{-1}^{1}\left[\sin(2ct)\right]_{-1-x}^{1-x} dx\\
&\geq &\frac{2 W^2}{3}-\frac{ W^2}{6c^2}(\sin(2c))^2
\end{eqnarray*}
On the other hand, from the proof of Lemma 2 of \cite{Bonami-Karoui2}, it can be easily checked that 
\begin{eqnarray}
 \|\mathcal {Q}_c\|^2_{HS} &\geq& \frac{2c}{\pi}-\frac{1}{\pi^2}\log\Big(\frac{2c}{\pi}\Big)-0.45\label{9}
\end{eqnarray}
By combining the previous two inequalities, one gets
$$
     \|\widetilde{Q}_{W,N}\|^2_{HS} \geq \frac{2c}{\pi}-\frac{1}{\pi^2}\log\big(\frac{2c}{\pi}\big)-0.45+\frac{2 W^2}{3}-
\frac{ W^2}{6c^2}(\sin(2c))^2.
$$
Since $Trace(\widetilde{Q}_{W,N})=2N W=\frac{2c}{\pi}$, then by using the previous inequality, one gets
\begin{eqnarray}
 Trace(\widetilde{Q}_{W,N}) - \|\widetilde{Q}_{W,N}\|^2_{HS} &=& \ds\sum_{k=0}^{N-1}\widetilde{\lambda}_{k,N}(W)(1-\widetilde{\lambda}_{k,N}(W))
\nonumber\\
&\leq& \frac{1}{\pi^2}\log\big(\frac{2c}{\pi}\big)+0.45-\frac{2 W^2}{3}+
\frac{ W^2}{6c^2}(\sin(2c))^2.
\end{eqnarray}
That is for $ c= \pi N W,$ we have 
\begin{equation}\label{10}
\eta(N,W)= \ds\sum_{k=0}^{N-1}\widetilde{\lambda}_{k,N}(W)(1-\widetilde{\lambda}_{k,N}(W))
\leq\frac{1}{\pi^2}\log(2 N W)+0.45-\frac{2 W^2}{3}+
\frac{ W^2}{6c^2}(\sin(2c))^2.
\end{equation}
Finally, since $\forall \; \varepsilon\in(0,\frac{1}{2})$  and $ x\in (\varepsilon,1-\varepsilon)$ , we have $x(1-x)\geq \varepsilon(1-\varepsilon)$ ,
then
\begin{eqnarray*}
  \varepsilon(1-\varepsilon) \mathcal N(W,\varepsilon) &\leq&\ds\sum_{k=0}^{N-1}\widetilde{\lambda}_{k,N}(W)(1-\widetilde{\lambda}_{k,N}(W))\leq \eta(N,W)
.\end{eqnarray*}
This conclude the proof of the theorem. $\qquad \Box$

\begin{remark}
We should mention that the upper bound  given by \eqref{Eigenvaluesnumber} outperforms the bound given in \cite{zhu}
in the sense that 
$$\mathcal N(W,\varepsilon) < \frac{\frac{2}{\pi^2}\log(N-1)+\frac{2}{\pi^2}\frac{2N-1}{N-1}}{\varepsilon(1-\varepsilon)},\quad 
\forall \; W\in \Big(0,\frac{1}{2}\Big),\quad N\geq 1.$$
\end{remark}

\begin{remark} We should mention that our estimate of  $\mathcal N(W,\varepsilon),$ the number of eigenvalues in the interval 
$(\varepsilon, 1-\varepsilon)$ and  given by \eqref{Eigenvaluesnumber}, is a non-asymptotic. It makes sense only if $\varepsilon $ is not too small.
 Recently, in \cite{Karnik}, the authors have given the following  asymptotic estimate of   $\mathcal N(W,\varepsilon),$ which is valid for small values of $\varepsilon,$ 
$$
\mathcal N(W,\varepsilon) = \left(\frac{8}{\pi^2} \log(8N+12)\right)\log\left(\frac{15}{\varepsilon}\right).
$$
\end{remark}

The following theorem is one of the main results of this work. It gives a fairly good bound of each eigenvalue $\widetilde{\lambda}_{n,N}(W)$ in terms  the corresponding eigenvalue $\lambda_n(c),$ with $c=\pi N W$ and $0\leq n\leq N-1.$ This allows us to generalize at ounce the various existing upper bounds for the  classical eigenvalues $\lambda_n(c).$

\begin{theorem}
Under the previous notation, for any integer $N\geq 1$ and real $W\in (0,1/2),$ we have for $c= N\pi W,$
\begin{equation}\label{Comparaison1}
\widetilde{\lambda}_{n,N}(W) \leq  A_W \,\,\,\lambda_n(c), \quad 0\leq n\leq N-1,
\end{equation}
where 
\begin{equation}\label{A}
 \frac{\pi^2}{8}\leq A_W = \frac{2 \pi^2}{\cos^2(\pi W)} \left(\frac{1}{4}-W^2\right)^2 \leq 2.
\end{equation}
\end{theorem}

\noindent
{\bf Proof: } We first use a classical technique for the construction of a subspace of the classical band-limited functions 
$$\mathcal B_{N\pi}=\{ f\in L^2(-\pi,\pi),\,\, \mbox{Supp}^t \widehat f \subseteq [-\pi,(2N-1)\pi]\}.$$ This is done as follows. Let $\varphi(\cdot)\in L^2(\mathbb R),$ with $\mbox{Supp}^t \widehat \varphi \subseteq [-\pi,\pi]$ and let 
$$ V_{N,\varphi}= \mbox{Span} \left\{ e^{2i \pi k t}\, \varphi(t),\,\, 0\leq k\leq N-1 \right\}.$$
That is if $f\in V_{N,\varphi},$ then ${\displaystyle f(t)= \sum_{k=0}^{N-1} \widehat P(k)  e^{2i\pi k t} \varphi(t) }.$ Here, 
${\displaystyle \sum_{k=0}^{N-1} \widehat P(k)  e^{2i\pi k t} = P(e^{2i\pi t}),}$ where $P\in \mathbb R_{N-1}[x]$ is a  polynomial of degree $N-1.$ Since 
$\mbox{Supp}^t \varphi \subseteq [-\pi,\pi]$ and since ${\displaystyle \widehat f(\xi)= \sum_{k=0}^N \widehat P(k) \widehat \varphi(\xi-2\pi k)},$ then 
$\mbox{Supp}^t \widehat f \subseteq [-\pi, (2N-1)\pi],$ that is $f\in \mathcal B_{N\pi}.$ By using Plancherel's equality, one gets 
$$\| f\|^2_{L^2(\mathbb R)} = \frac{1}{2\pi} \|\widehat f\|^2_{L^2(\mathbb R)}=\frac{1}{2\pi} \sum_{k=0}^{N-1} |\widehat P(k)|^2 \|\widehat \varphi\|^2_{L^2(\mathbb R)}.$$
Also, from Parseval's equality, we have
$$ \sum_{k=0}^{N-1} |\widehat P(k)|^2 = \| P(e^{2i\pi t})\|^2_{L^2(-1/2,1/2)}.$$
By combining the previous two equalities, one gets
$$ \| f\|^2_{L^2(\mathbb R)}=\frac{1}{2\pi} \| P(e^{2i\pi t})\|^2_{L^2(-1/2,1/2)} \|\widehat \varphi\|^2_{L^2(\mathbb R)},\quad \deg P\leq N-1.$$
On the other hand, for $W\in (0, 1/2),$ we have
\begin{eqnarray*}
 \| f\|^2_{L^2(-W,W)}&=& \int_{-W}^W  | P(e^{2i\pi t})|^2 \varphi^2(t) \, dt \geq \min_{t\in [-W,W]} |\varphi(t)|^2 
 \int_{-W}^W  | P(e^{2i\pi t})|^2\, dt.
\end{eqnarray*}
Hence, for any $f\in V_{N,\varphi},$ we have 
\begin{equation}\label{Ineq1}
\frac{1}{\min_{t\in [-W,W]}|\varphi^2(t)|} \frac{\| f\|^2_{L^2(-W,W)}}{ \| f\|^2_{L^2(\mathbb R)}}\geq 2\pi \frac{\| P(e^{2i\pi t})\|^2_{L^2(-W,W)}}{\| P(e^{2i\pi t})\|^2_{L^2(-1/2,1/2)} \|\widehat \varphi\|^2_{L^2(\mathbb R)}}.
\end{equation}
In particular, for ${\displaystyle \widehat \varphi(\xi)=\mathbf 1_{[-\pi,\pi]}(\xi) \cos(\xi/2)},$ we have ${\displaystyle \|\widehat \varphi\|^2_{L^2(\mathbb R)}=\pi.}$ Moreover, we have 
$$\varphi(x)= \frac{1}{2\pi} \int_{-\pi}^\pi e^{i x\xi} \cos(\xi/2) \, d\xi= \left\{\begin{array}{ll} \frac{1}{2\pi} \left(\frac{\cos(\pi x)}{\frac{1}{4} - x^2}\right), & x\in (-1/2,1/2)\\ & \\ \frac{1}{2} & \mbox{if } x=\pm \frac{1}{2}.\end{array}\right.$$
So that 
$$ \min_{x\in [-W,W]} |\varphi(x)|^2= \frac{1}{4 \pi^2} \left(\frac{\cos(W \pi)}{\frac{1}{4}-W^2}\right)^2.$$
Hence, for this choice of $\varphi$ and by using \eqref{Ineq1}, one concludes that for any  polynomial $P_N\in \mathbb R_{N-1}[x]$ of
 degree $N-1,$ we have
\begin{equation}\label{Ineq2}
\frac{\| P_N(e^{2i\pi t})\|^2_{L^2(-W,W)}}{\| P_N(e^{2i\pi t})\|^2_{L^2(-1/2,1/2)}} \leq 2 \left(\frac{\frac{1}{4}-W^2}{\cos(W\pi)}\right)^2 \frac{\|f_N\|^2_{L^2(-W,W)}}{\|f_N\|^2_{L^2(\mathbb R)}},
\end{equation}
where, ${\displaystyle f_N(t)= P_N(e^{2i\pi t}) \varphi(t)\in V_{N,\varphi}.}$ Next, let $S_N$ be the subspace of sequences 
 ${\pmb x}=(x_n)_n \in l^2(\mathbb C)$ with elements indexed on 
$[[0,N-1]],$ so that ${\displaystyle  \widehat{{\pmb x}}(t)= \sum_{k=0}^{N-1} x_k e^{2i\pi k t}.}$ Also, we denote by $s_n, v_n,$  the $(n+1)-$dimensional subspace of $S_N$ and $V_{N,\varphi},$ respectively. Note that  the eigenvalues of the Sinc-kernel operator, are invariant under dilation of the time-concentration interval and translation and dilation of the bandwidth concentration interval. That is  for $\tau, c>0,$ we have 
$\lambda (\mathcal Q_{\tau,c})=\lambda(\mathcal Q_{1,\tau c})=\lambda(\mathcal Q_{\tau c})$ or equivalently 
$\lambda_{n,\tau}(c)= \lambda_{n,1}(\tau c)=\lambda_n( \tau c).$ By using the previous properties as well as the Min-Max characterisation of this later, together with inequality \eqref{Ineq2}, and the fact that $V_{N,\varphi}$ is a subspace of $\mathcal B_{(N+1)\pi},$ one gets 
\begin{eqnarray*}
\widetilde{\lambda}_{n,N}(W)=\max_{S_n}\min_{{\pmb x}\in S_n\setminus \{0\}}\frac{\|\widehat{{\pmb x}}\|^2_{L^2(-W,W)}}{\| {\widehat {\pmb x}}\|^2_{L^2(-1,1)}}&\leq & A_{W}
\max_{v_n}\min_{f\in v_n\setminus \{0\}} \frac{\|f\|^2_{L^2(-W,W)}}{\|f\|^2_{L^2(\mathbb R)}}\\
&\leq& A_{W} \max_{U_n}\min_{f\in U_n\setminus \{0\}} \frac{\|f\|^2_{L^2(-W,W)}}{\|f\|^2_{L^2(\mathbb R)}}\\
&\leq & A_{W} \max_{W_n}\min_{f\in W_n\setminus \{0\}} \frac{\|f\|^2_{L^2(-1,1)}}{\|f\|^2_{L^2(\mathbb R)}}= A_{W} \lambda_n(c).
\end{eqnarray*}
Here, ${\displaystyle  A_W = \frac{2 \pi^2}{\cos^2(\pi W)} \left(\frac{1}{4}-W^2\right)^2},$   $U_n$ is an $(n+1)-$dimensional subspace of $\mathcal B_{N\pi}$ and $W_n$ is an $(n+1)-$subspace of 
${\displaystyle \mathcal B_c= \{f\in L^2(\mathbb R),\, \mbox{Supp}^t \widehat f \subseteq [-c,c]\},\,\, c= \pi N W.}$ $\qquad \Box$\\

The previous theorem allows us to extend some known estimates for the classical $\lambda_n(c)$ to the eigenvalues 
$\widetilde \lambda_{n,N}(W).$ In  \cite{BJK}, it has been shown that for any $c>0$ and for any $n\geq \max\Big(2,\frac{ec}{2}\Big),$ we have 
${\displaystyle \lambda_n(c) \leq e^{-(2n+1)\log(\frac{2}{ec}(n+1))}}.$  By using the previous theorem, together with the previous  non-asymptotic estimate of the $\lambda_n(c),$  one concludes that for any $N\geq 3$ and 
any $W\in \big(0,\frac{2}{e\pi}\frac{N-1}{N}\big),$ we have
	\begin{equation}\label{Eq}
	\widetilde{\lambda}_{k,N}(W)\leq 2 e^{-(2k+1)\log(\frac{2}{e\pi NW}(k+1))},\qquad 2\leq \frac{e\pi}{2} NW\leq k\leq N-1,	
	\end{equation}
for any $N\geq 3$ and 	 $W\in \big(0,\frac{2}{e\pi}\frac{N-1}{N}\big).$ Moreover, it has been shown in \cite{BJK} that  for any $\frac{2c}{\pi}+\log(c)+6\leq n\leq c,$ there exists a uniform constant $\eta >0$ such that 
${\displaystyle \lambda_n(c)\leq \exp\left(-\eta \frac{n-\frac{2c}{\pi}}{\log(c)+5}\right).}$ This last estimate combined with the previous theorem, give us the following similar estimate
\begin{equation}\label{ddecay7}
\widetilde{\lambda}_{n,N}(W) \leq  2 e^{- \eta \frac{n-2NW}{\log(\pi NW)+5}},\qquad 2NW +\log(\pi N W)+6\leq n\leq \pi N W.
\end{equation}
It is interesting to note that besides providing an explicit exponential decay rate for the $\widetilde \lambda_{n,N}(W),$ the estimate \eqref{Eq} 
provides us with  estimates for the unknown constants $C_1(W,\varepsilon),$ $C_2(W,\varepsilon)$  appearing in the 
following asymptotic decay rate, given in \cite{Slepian}
\begin{equation}\label{Decay1}
\widetilde \lambda_{k,N}(W) \leq C_1(W,\varepsilon) e^{-C_2(W,\varepsilon) N},\quad \forall\, k\geq \lceil{2NW(1+\varepsilon)}\rceil,\quad N \geq N_1(W,\varepsilon).
\end{equation}
More precisely, by comparing \eqref{Eq} and \eqref{Decay1}, one concludes that for  $N\geq 3,$  $\varepsilon >\frac{e\pi -6}{4}$ and
$W \leq \frac{2}{e\pi}\frac{N-1}{N},$ we have 
\begin{equation}\label{estimates}
C_1(W,\varepsilon) \leq 2,\qquad C_2(W,\varepsilon) \geq 4 W(1+\varepsilon)\log\left(\frac{4(1+\varepsilon)+2}{e\pi}\right).
\end{equation}

\section{Applications.}

In this paragraph, we give two applications of the DPSWF's. The first application is related to  a lower bound  estimate for the constant appearing in 
the  Tur\`an-Nazarov concentration  inequality, see \cite{Nazarov}. The second applications deals with the quality of approximation by the DPSWFs of Bandlimited functions and functions from periodic Sobolev spaces. \\

Let us first recall the following Tur\`an-Nazarov type concentration  inequality. Let $\TT$ be the unit circle and let $\mu$ be the Lebesgue measure on $\TT,$ normalized so that $\mu(\TT)=1,$ then for every $0\leq q \leq 2,$ every 
trigonometric polynomial 
$$
{\displaystyle P(z)=\sum_{k=1}^{n+1} a_k z^{\alpha_k},\quad a_k\in \mathbb C,\quad \alpha_k \in \mathbb N,\quad  z\in \TT,}
$$
and every measurable subset $E\subset \TT,$ with $\mu(E)\geq \frac{1}{3},$ we have
\begin{equation}
\label{Ineq_Turan}
\| P \|_{L^q(E)} \geq e^{-A \, n\,  \mu(\TT\setminus E)}  \| P\|_{L^q(\TT)}.
\end{equation}
Here, $A$ is a constant independent of $q,$ $E$ and $n.$ Since, the DPSWFs $U_{n,W}^N$ are given by 
 $$U_{n,W}^N(x)=\epsilon_n\sum_{k=0}^{N-1}v_k^n(e^{-2i\pi x})^{\frac{N-1}{2}-k}=P_N(e^{-2i\pi x}),\quad x\in [-1/2,1/2],$$
then by combining and inequality \eqref{Ineq_Turan} with  $q=2$ and $E=(-W,W)$ where $1/6\leq W < 1/2$, one gets	
	\begin{equation} 
\widetilde \lambda_{n,N}(W)=\frac{\norm{U_{n,W}^N}_{L^2(-W,W)}}{\norm{U_{n,W}^N}_{L^2(-1/2,1/2)}}\geq e^{-A(1-2W)(N-1)},\quad \forall\, 0\leq n\leq N-1.
	\end{equation}
In particular, for $n=N-1,$  $W=\frac{1}{6}$ and by using the estimate \eqref{Eq}, together with a straightforward computation, one gets 
\begin{equation}\label{estimateA}
A \geq \frac{2}{1-2/6}\log\left(\frac{12}{e \pi}\right) = 1.02.
\end{equation} 

Concerning the quality of approximation of bandlimited functions by the DPSWF's,  we have a partial result. 
In fact,  we check that under some conditions on $W$ and $N,$ there exists $N_0\geq [2NW]-1,$ such that  the eigenspace spanned by the first $N_0$ dilated DPSWFs $\sqrt{W} U_{k,W}^N(W \cdot)$ approximates  the eigenspace spanned by the corresponding classical $\psi_{k,c}.$ For this purpose, we 
first recall the following result given by Theorem 3 of \cite{Zwald} and concerning the approximation of eigenspaces spanned by a set of eigenfunctions
of positive self-adjoint Hilbert-Schmidt operator and its positive  self-adjoint perturbed version. More precisely, if $A$ is such an operator with simple eigenvalues $\lambda_0>\lambda_1>\cdots $ and if there exists an integer $D>0$ such that $\lambda_D>0$ and $\delta_D=\frac{1}{2}(\lambda_D-\lambda_{D+1})$ and if $A+B$ is such a perturbed operator satisfying the extra condition that $\| B\|<\delta_D/2,$ then 
\begin{equation}\label{projections}
\| \pi_D(A) - \pi_D(A+B)\| \leq \frac{\|B\|}{\delta_D}.
\end{equation} 
Here, $\pi_D(A)$ denotes the orthogonal projection over the space spanned by the first eigenfunctions of the operator $A.$ 
In the sequel, we let $\widetilde L_{W,N}$ denote the operator defined on $L^2(-1,1)$ by
\begin{equation}\label{OperatorL}
\widetilde L_{W,N}(f)(x)=\int_{-1}^1 \frac{\sin(\pi N W(x-y))}{\sin(\pi W (x-y))} f(y)\, dy,\quad x\in (-1,1).
\end{equation}
Then it is easy to check that the $N$ dilated DPSWF's $\sqrt{W} U_{k,W}^N(W \cdot)$ are eigenfunctions of $\widetilde L_{W,N},$ with the same associated 
eigenvalues $\widetilde \lambda_{k,N}(W)$ as the usual DPSWF's. 
In the special case where 
for $c=\pi NW,$ the operators $A$ and $A+B$ are  given by $\mathcal Q_c,$  and $\widetilde L_{W,N},$ respectively, we obtain the following proposition that gives us an approximation of eigenspaces spanned by classical PSWFs and the corresponding DPSWFs.

\begin{proposition}\label{ApproxEigen}
 Let $\pi_K$ and $\widetilde \pi_K$ be the two projection operators on the spaces spanned by the first $K$-eigenfunctions of the the  operator $\mathcal Q_{c}$ and  $\widetilde{L}_{W,N},$ respectively.
For any real $b>\frac{\log 3}{\pi},$ there exists  $c_b >0$ such that for any $(W, N)\in \big(0,\frac{1}{2}\big)\times \mathbb N,$ with 
\begin{equation}\label{conditionNW}
c_b \leq \pi N W \leq \exp\left(\alpha_b \frac{\sin 2\pi W}{W^3}-2\log 2 -\frac{\pi}{b}\right),\qquad \alpha_b=\frac{3}{32 b\pi}\Big(1-\frac{3}{1+e^{\pi b}}\Big),
\end{equation}
then there exists $N_0\geq [2NW]$ such that 
\begin{equation}\label{Eigenspaces}
\|\pi_{N_0}-\widetilde \pi_{N_0}\|\leq W^3 \left(\frac{4b\pi}{3 \sin(2\pi W)}\right)\frac{\log(\pi N W)+2\log 2 +\frac{\pi}{b}}{1-\frac{3}{1+e^{\pi b}}}.
\end{equation} 
\end{proposition}

\noindent
{\bf Proof: } 
We first recall that in \cite{Slepian2} and for a fixed $b\geq 0,$  $c>0,$   the author has  given the following limit  result for  $\lambda_{n}(c),$
\begin{equation}\label{limit1}
 \lim_{c\rightarrow +\infty} \lambda_{n_{c,b}} (c) = \frac{1}{1+e^{\pi b }},\qquad n_{c,b} = \left[ \frac{2c}{\pi}+\frac{2 b}{\pi} \log 2+\frac{b}{\pi} \log c\right]
\end{equation} 
Hence, by applying the previous estimate for the two fixed values of $b=0$ and $b> \frac{\log 3}{\pi},$ one concludes that there exists $C_b>0$ such that for any $c\geq c_b,$ we have 
\begin{eqnarray}\label{Estimate1}
0<\sum_{k=n_{c,0}}^{n_{c,b}-1} \lambda_{k}(c)- \lambda_{k+1}(c)&=&\lambda_{n_{c,0}}(c)-\lambda_{n_{c,b}} (c) \leq \frac{1}{2}-\frac{3}{2} \frac{1}{1+e^{\pi b}}.
\end{eqnarray}
Note that from  \eqref{limit1}, we have 
$$ n_{c,b}-n_{c,0}=\left[ \frac{2c}{\pi}+\frac{2 b}{\pi} \log 2+\frac{b}{\pi} \log c\right]-\left[ \frac{2c}{\pi}\right]\leq \frac{b}{\pi} \log c+\frac{2 b}{\pi} \log 2+1.$$
Consequently, by using \eqref{Estimate1}, one concludes  there exists $N_0\geq n_{c,0}$ such that 
$$\delta_{N0}= \frac{\lambda_{N_0}(c)- \lambda_{N_0+1}(c)}{2}\geq \frac{1}{\frac{b}{\pi} \log c+\frac{2 b}{\pi} \log 2+1}\left(\frac{1}{4}-\frac{3}{4} \frac{1}{1+e^{\pi b}}\right).$$
Next, we consider  the special cases of $c= \pi N W,$ and the operators $A$ and $A+B$ are given by $\mathcal Q_c,$ $\widetilde L_{W,N},$ respectively.
By using  \eqref{normsHS} and  \eqref{conditionNW}, one can easily check that  
$$ \| B\|= \|\widetilde L_{W,N}- \mathcal Q_c\|\leq \|\widetilde L_{W,N}- \mathcal Q_c\|_{HS}\leq \frac{\delta_{N0}}{2}.$$
Hence by using \eqref{projections} and \eqref{normsHS}, one gets the desired result \eqref{Eigenspaces}.$\qquad \Box$

It is well known  see for example \cite{Sl1, Wang1}, that if $f\in B_c,$ where $B_c$ is the space of bandlimited functions, given by 
\eqref{Bc},  then we have 
 \begin{equation}\label{eq6.0}
 \|f- \pi_{N_0} f\|_{L^2(-1,1)} \leq   \lambda_{N_0}(c) \|f\|_{L^2(\mathbb R)}.
\end{equation}
Here, $\pi_{N_0}$ is the orthogonal projection over the first classical PSWFs $\psi_{k,c}(\cdot).$
\begin{remark}\label{EigenspacesApprox}
By combining \eqref{Eigenspaces} and the previous inequality, one gets the following partial result concerning the quality of approximation 
of bandlimited functions by the dilated  DPSWF's. For $c=\pi N W,$   $b>\frac{\log 3}{\pi}$ and under condition 
\eqref{conditionNW}, there exists $N_0 \geq [2NW]$ such that for $f\in B_c,$ we have 
 \begin{equation}\label{eq6.1}
 \|f- \widetilde \pi_{N_0} f\|_{L^2(-1,1)} \leq  W^3 \left(\frac{4b\pi}{3 \sin(2\pi W)}\right)\frac{\log(\pi N W)+2\log 2 +\frac{\pi}{b}}{1-\frac{3}{1+e^{\pi b}}} \|f\|_{L^2(-1,1)}+ \lambda_{N_0}(c) \|f\|_{L^2(\mathbb R)}.
\end{equation}
Here, $\widetilde \pi_{N_0}$ is the orthogonal projection over the first $N_0$ dilated DPSWFs $\sqrt{W} U^N_{k,W}(W \cdot).$ In a similar manner, we may 
extend this approximation quality of the dilated DPSWF's in the more general class of functions  of almost time- and band-limited functions. For more details on this class of functions, the reader is refereed to \cite{JKS, Landau}.  We leave the details of this extension to the reader.
\end{remark}

Next,  check that the DPSWFs are well adapted for the spectral approximation of functions from the periodic Sobolev 
spaces.  Note that for a given real number $s>0,$  the periodic 
Sobolev space $H_{per}^s{([-1/2,1/2])}$ is defined by
$$ H_{per}^s{(-1/2,1/2)}=\left\{ f\in {L^2(-1/2,1/2)},\, \|f\|_{H^s}^2 =\sum_{n \in \mathbb{Z}}(1+n^2)^s|\hat{f}_n|^2 <+\infty \right\},$$ where
${\displaystyle \widehat f_n= \int_{-1/2}^{1/2} f(x) e^{-2i\pi nx}\, dx.}$

\begin{lemma}
let $W\in (0,\frac{1}{2}),$ $N\in \mathbb N$ such that $\pi N W \geq 1.$ Let  $s>0,$  and  $H_{per}^s{(-1/2,1/2)},$ then there exists a constant $M>0$ such that for any integer  $  [2NW]+\log(\pi NW) +6\leq K \leq N-1,$ we have 
\begin{equation}\label{approx}
    \|f-\widetilde \pi_K(f)\|_{L^2(-W,W)}\leq\frac{4}{(4+N^2)^{s/2}}\|f\|_{H^s}
+ \sqrt{\widetilde \lambda_{K,N}(W)} \|f\|_{L^2(-1/2,1/2)}.
\end{equation}
Here, $\widetilde \pi_K$ is the orthogonal projection over the space spanned by the first $K$ DPSWFs, associated with the parameters $W, N.$
\end{lemma}

\noindent
{\bf Proof:} We first note that if $f_N\in L^2(-1/2,1/2)$ is the function given by ${\displaystyle f_N(x) = \sum_{k=-[(N-1)/2]}^{[(N-1)/2]} \widehat f_k e^{2i\pi k x},}$ then we have 

\begin{equation}\label{Eq1}
 \|f-f_N\|^2_{L^2(-1/2,1/2)}= \ds\sum_{|n|\geq [(N+1)/2]}^{+\infty}\frac{1}{(1+n^2)^s}(1+n^2)^s|\hat{f}_n|^2\leq \frac{4}{(4+N^2)^s}
\|f\|^2_{H^s}.
\end{equation}
On the other hand, by using the expressions of the DPSWFs, given by \eqref{DPSWFs}, as well as their double orthogonality property \eqref{doubleorthogonality}, on gets 
\begin{equation}\label{12}
   f_N(x)=\ds\sum_{k=0}^{N-1}\beta_k U_{k,W}^N(x) ,\;\;\;\forall x \in   [-1/2,1/2],\qquad 
  f_N(x) =\ds\sum_{k=0}^{N-1} \alpha_k \frac{U_{k,W}^N(x)}{\sqrt{\widetilde\lambda_{k,N}(W)}} ,\;\;\;\forall x \in   [-W,W].
\end{equation}
Since the previous two expansions coincide on  $[-W,W],$ then we have 
\begin{eqnarray}\label{14}
 \alpha_k=\beta_k\sqrt{\widetilde \lambda_{k,N}(W)}
\end{eqnarray}
Moreover by using the previous identity, together with  Parseval's equality and the decay of the $\widetilde \lambda_{k,N}(W),$ one gets
\begin{eqnarray}\label{16}
\|f_N-\widetilde \pi_K(f_N)\|^2_{L^2(-W,W)}&=&\ds\sum_{k=K}^{N-1}|\alpha_{k}|^2 = \ds\sum_{k=K}^{N-1} \widetilde\lambda_{k,N}(W) |\beta_k|^2 \nonumber \\
&\leq & \widetilde \lambda_{K,N}(W) \ds\sum_{k=K}^{N-1} |\beta_k|^2 \leq \widetilde \lambda_{K,N}(W)|\|f\|^2_{L^2(-1/2,1/2)}.
\end{eqnarray}
Moreover, since $\widetilde \pi_K$ is an orthogonal projection, then we have $\|\widetilde \pi_K\|\leq 1.$ Hence, by using \eqref{Eq1} and 
\eqref{16}, one gets 
\begin{eqnarray}    % \nonumber to remove numbering (before each equation)
       \|f-\widetilde \pi_K(f)\|_{L^2(-W,W)}&\leq&  \|f-f_N\|_{L^2(-W,W)}+ \|\widetilde \pi_K (f-f_N)\|_{L^2(-W,W)}+\|f_N-\widetilde \pi_K(f_N)\|_{L^2(-W,W)} \nonumber\\
 &\leq&  2 \|f-f_N\|_{L^2(-W,W)}  +\sqrt{\widetilde \lambda_{k,N}(W)}  \|f\|^2_{L^2(-1/2,1/2)} \nonumber\\
 &\leq&  \frac{4}{(4+N^2)^{s/2}} \|f\|_{H^s}  +\sqrt{\widetilde \lambda_{k,N}(W)}  \|f\|^2_{L^2(-1/2,1/2)}.   \qquad \Box 
\end{eqnarray}

\begin{remark}
It is easy to check that by considering the dilated DPSWF's $\sqrt{W} U^N_{k,W}(W\cdot)$ and by considering the  periodic extension of $f\in H^s(-1,1),\, s>0,$ we get the following  approximation result of $f$ by the first $K$ dilated DPSWF's,
\begin{equation}\label{approx2}
    \|f-\widetilde \pi_K(f)\|_{L^2(-1,1)}\leq\frac{4}{(4+N^2)^{\frac{s}{2}}}\|f\|_{H^s}
+ \sqrt{\widetilde \lambda_{K,N}(W)} \|f\|_{L^2\big(-1/(2W),1/(2W)\big)}.
\end{equation}
\end{remark}

\section{Numerical results.}

In this section,  we give three examples that illustrate the different results of this work. \\

\noindent
{\bf Example 1:} In  this first example, we  give different  numerical tests that illustrate the result 
of Proposition~1, implying in particular that  each eigenvalue $\widetilde \lambda_{k,N}(W)$ is well approximated by the corresponding 
classical $\lambda_k(c),\,\, c=\pi N W.$ Also, these tests illustrate the unexpected and important 
inequality \eqref{Comparaison1} of Theorem 2 that bounds each  $\widetilde \lambda_{n,N}(W)$ in terms of the corresponding  $\lambda_n(c),\,\, c=\pi NW$
and up to a small constant $A_W.$ This  
allows us also to check the exponential decay rates for the $\widetilde \lambda_{n,N}(W),$   given by \eqref{ddecay7} when  $n$ is close to the plunge region around $[2NW],$  and by \eqref{Eq} when $n$ is sufficiently far from $[2NW]$. For these purposes,
we have considered the value of $N=60$ and the four values of $W=0.1,\, 0.2,\, 0.3,\, 0.4.$  Note that  the  $\widetilde \lambda_{n,N}(W)$ 
are computed  by computing the eigenvalues of the Toeplitz matrix \eqref{Toeplitzmatrix}. The corresponding eigenvalues $\lambda_n(c)$ associated with the classical PSWFs are computed with high precision  by using the method given in \cite{Karoui-Moumni}. In Table 1, we have listed the $\ell_2$-approximation error 
corresponding to both sequences of eigenvalues as predicted by proposition 1. In Figure 1(a), we have plotted the graphs of $\widetilde \lambda_{n,N}(W)$ 
for the considered four values of $W.$ This figure illustrate Theorem 1 in the sense that the sequence $\widetilde \lambda_{n,N}(W)$ clusters around $1$ and $0$ and the number of the eigenvalues in the plunge region of the spectrum follows the bound given by Theorem 1.  Finally,  to illustrate the exponential decay rate of the $\widetilde \lambda_{n,N}(W)$ as well as the main result of Theorem 2, we have plotted in Figure 1(b) the graphs of $\log\big(\widetilde \lambda_{n,N}(W)\big)$ versus the corresponding $\log\big(\lambda_n(c)\big),\, c=\pi N W.$  

\begin{center}
  \begin{table}[h]
  \vskip 0.2cm\hspace*{3cm}
  \begin{tabular}{ccc} \hline
   $W$ &$c= \pi N W$&$\| \lambda(\mathcal Q_c)-\lambda(\widetilde Q_{W,N})\|_{\ell_2}$\\   \hline
  $0.1$  &$18.85$   &$4.15E-03 $  \\
  $0.2$ &$37.70 $   &$1.65E-02 $   \\
  $0.3$ &$56.55 $   &$3.98E-02 $   \\
  $0.4$ &$75.40 $   &$8.51E-02 $   \\ \hline
     & &  \\
  \end{tabular}
  \caption{Illustration of Proposition~1 concerning the $\ell_2$-error approximation of the 
  sequence $(\widetilde \lambda_{n,N}(W))_n$ by the sequence $(\lambda_n(c))_n$ for $N=60$ and different values of $W.$ }
  \end{table}
  \end{center}

 \begin{figure}[h]\hspace*{0.05cm}
 {\includegraphics[width=15.05cm,height=5.5cm]{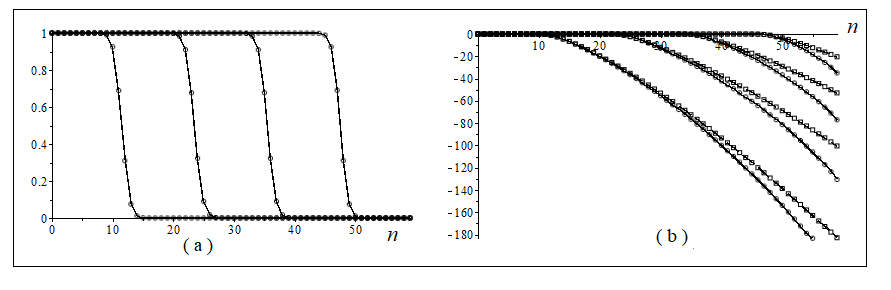}}
 \vskip -0.5cm\hspace*{1cm} \caption{(a) Graphs of the eigenvalues $\widetilde \lambda_{n,N}(W)$ for the  values of $W=0.1,\, 0.2,\, 0.3,\, 0.4$ (from left to right) and  different values of $n.$ 
 (b) Graphs of the associated $\log\big(\widetilde \lambda_{n,N}(W)\big)$ (in circles) versus the corresponding $\log\big(\lambda_n(c)\big),\, c=\pi N W.$ (in boxes) }
 \end{figure}

\noindent
{\bf Example 2:} In this second example, we illustrate the quality of the spectral  approximation of bandlimited functions by the DPSWfs, as partially predicted by Proposition~\ref{ApproxEigen} and Remark~\ref{EigenspacesApprox}. For this purpose, we have considered the $\alpha$-bandlimited function $f_\alpha$ defined by 
${\displaystyle f_\alpha(x)=\frac{\sin(\alpha x)}{\alpha x} }$ with $\alpha =56$ and  the special values of $W=0.3$ and $N=60,$ so that $c_N=\pi N W=56.55.$ Then,  we have computed $\widetilde \pi_N f_\alpha,$ the orthogonal projection of $f_\alpha$  over the finite dimensional subspace spanned by the orthonormal set 
of $L^2(-1,1)$ given by the dilated DPSWFs, $(\sqrt{W} U^N_{k,W}(W \cdot))_{0\leq k\leq N-1}.$ We found that
$$\sup_{t\in [-1,1]} |f_\alpha(t)-\widetilde \pi_N f_\alpha(t)| \approx 4 E-10.$$ That is  $\widetilde \pi_N f_\alpha$ provides us with  a surprising high approximation of the  $f_\alpha.$  \\

\noindent
{\bf Example 3:} In this last example, we illustrate the quality of approximation of function from the periodic Sobolev space 
$H^s_{per}(-1,1)$ by the dilated DPSWFs, as predicted by Lemma~2 and Remark~3. For this purpose,  we consider   the
Weierstrass function
\begin{equation}\label{eq5.1}
\mathcal W_s(x)= \sum_{k\geq 0} \frac{\cos(2^{k}x)}{2^{ks}},\quad -1\leq
x\leq 1.
\end{equation}
Note that $\mathcal W_s \in H_{per}^{s-\epsilon}(-1,1),\,\forall 0<\epsilon <s.$ We consider the special value of $s=1$ and the  same set of 
dilated DPSWFs $(\sqrt{W} U^N_{k,W}(W \cdot))_{0\leq k\leq N-1},$ of the previous example with $W=0.3$ and $N=60.$ Then, we have computed
the orthogonal projections $\widetilde \pi_K \mathcal W_s,$ over the subspace spanned by the first $K$ dilated DPSWFs. We found that 
$$\| \mathcal W_1 - \widetilde \pi_{K} \mathcal W_1\|_{L^2(-1,1)}\approx 8.64 E-03,\quad K=N=60$$
an
$$\| \mathcal W_1 - \widetilde \pi_{K} \mathcal W_1\|_{L^2(-1,1)}\|\approx 2.43E-02,\quad K=[2NW]=36.$$
Note that this Weierstrass function has been already used in \cite{Bonami-Karoui3} to test the quality of approximation of functions from the Sobolev spaces $H^s(-1,1)$ by the classical PSWFs $\psi_{n,c}.$ The previous two approximation errors and the numerical results given in \cite{Bonami-Karoui3}, 
indicate that the DPSWFs outperform the classical PSWFs for this kind of spectral approximation. In fact, with  fewer expansion coefficients, for this example, the  DPSWFs provide better approximation of this  Weierstrass function. This indicates that a DPSWF's based scheme for the approximation of 
Sobolev spaces over compact intervals can be complementary to the proposed similar schemes based on the classical PSWFs, that have been studied in 
\cite{Bonami-Karoui3, Wang1, Wang2}.


\begin{thebibliography}{99}

\bibitem{Adcock} B. Adcock and D. Huybrechs, On the numerical stability of Fourier extensions,
{\it Foundations of Computational Mathematics,} {\bf 14} (4) (2014), 635--687.

\bibitem{Batir} N. Batir,  Inequalities for the gamma function. {\it  Arch. Math.} 2008; { 91}: 554--563.

\bibitem{Bonami-Karoui2} A. Bonami and A. Karoui, Random Discretization of the Finite Fourier Transform and Related Kernel Random Matrices, (2019)
availbale at https://arxiv.org/abs/1703.10459.

\bibitem{Bonami-Karoui1} A. Bonami and A. Karoui, Spectral Decay of Time and Frequency Limiting Operator, {\it Appl. Comput. Harmon. Anal.} 	{\bf 42} (2017), 1--20.

\bibitem{Bonami-Karoui3}  A. Bonami and A. Karoui, Approximation in Sobolev spaces by Prolate Spheroidal Wave Functions,
 {\it Appl. Comput . Harmon. Anal.,}  {\bf 42} (2017), 361--377.

\bibitem{BJK} A. Bonami, P. Jaming  and A. Karoui, Non-Asymptotic Behaviour of the Sinc-Kernel Operator and Related Applications,
 available at arXiv:1804.01257, (2018).
 
\bibitem{Wakin1} M. A. Davenport and M. B. Wakin, Compressive sensing of analog signals using Discrete Prolate Spheroidal Sequences,
{\it Appl. Comput. Harmon. Anal.} 	{\bf 33} (2012), 438--472.

\bibitem{JKS}  P. Jaming, A. Karoui, S. Spektor, The approximation of almost time- and band-limited
functions by their expansion in some orthogonal polynomials bases, {\it J. Approx. Theory,} {\bf 212} (2016), 41--65.


 \bibitem{HL}  J. A. Hogan and J. D. Lakey, {\it Duration and Bandwidth Limiting: Prolate Functions, Sampling, and Applications,}
 Applied and Numerical Harmonic Analysis Series, Birkh\"auser, Springer, New York, London, 2013.
 
 \bibitem{Landau} H. J. Landau and H. O. Pollak, Prolate spheroidal  wave functions, Fourier analysis and
 uncertainty-III. The dimension of space of essentially time-and  band-limited signals, {\it Bell System Tech. J.} {\bf 41}, (1962),
 1295--1336.
 
 \bibitem{Karnik} S. Karnik, Z. Zhu, M. Wakin, J. Romberg, and M. Davenport, The  Fast Slepian Transform,  {\it Appl. Comput . Harmon. Anal.,}  {\bf 46} (2019), 624--652. 
 
 \bibitem{Karoui-Moumni} A. Karoui and T. Moumni,  New efficient  methods of computing the prolate spheroidal wave functions and
 their corresponding eigenvalues, {\it Appl. Comput. Harmon. Anal.,}  {\bf 24}, No.3,  (2008), 269--289.

\bibitem{Karoui-Souabni} A. Karoui and  A. Souabni, Generalized Prolate Spheroidal Wave Functions: Spectral Analysis  and Approximation of Almost Band-limited Functions, {\it J.  Fourier  Anal.  Appl.}  {\bf 22} (2016), 383--412.

\bibitem{Nazarov}  F.\,L. Nazarov,
Complete Version of Tur\`an's Lemma for Trigonometric Polynomials on the Unit Circumference, in 
{\it Operator Theory: Advances and Applications,} {\bf 113}, Birkh\"aseur Verlag Basel, Switzerland (2000),  239--246.


\bibitem{NIST} F.W. Olver, D.W. Lozier, R.F.  Boisvert,
C.W. Clark, { NIST Handbook of Mathematical Functions,}  Cambridge University Press; New York;  2010.

\bibitem{Sl1}D. Slepian, H. O. Pollak, Prolate spheroidal wave functions, Fourier analysis and uncertainty I, Bell System Tech. J. 40 (1961), 43-64.

\bibitem{Slepian} D. Slepian,  Prolate spheroidal wave functions, Fourier analysis and
uncertainty--V: The Discrete Case, {\it Bell System Tech. J.} {\bf 57}
(1978), 1371--1430.

\bibitem{Slepian2} D. Slepian, Some Asymptotic Expansions for Prolate Spheroidal Wave Functions, {\it Stud. Appl. Math., } {\bf 44} (4), (1965), 99--140.

\bibitem{Wang1} L. L. Wang,  Analysis of spectral approximations using prolate spheroidal wave functions.
 Math. Comp. {\bf 79} (2010), no. 270, 807--827.

\bibitem{Wang2} L. L. Wang and J. Zhang, A new generalization of the PSWFs with applications to spectral
 approximations on quasi-uniform grids, {\it Appl. Comput. Harmon. Anal.}
 {\bf 29},  (2010), 303--329.


\bibitem{Feng} F. Yin, C. Debes and A. M. Zoubir, Parametric waveform design using Discrete Prolate Spheroidal Sequences for enhanced detection of extended targets, {\it IEEE Trans. Signal Process., } {\bf 60} (9) (2012), 4525--4536.

\bibitem{zhu} Z. Zhu, M. B. Wakin, Approximating Sampled Sinusoids and Multiband Signals Using Multiband Modulated DPSS Dictionaries, 
{\it J.  Fourier  Anal.  Appl.}  {\bf 23} (2017), 1263--1310.

\bibitem{Zwald} L. Zwald Laurent and G. Blanchard, On the Convergence of Eigenspaces in Kernel Principal Component Analysis, in Advances in Neural Information Processing Systems 18: Proceedings of the 2005 Conference (Neural Information Processing),  MIT Press, (2006).


\end{thebibliography}
\end{document}